\documentstyle[12pt]{article}
\begin{document}
\begin{titlepage}
\begin{flushright}
math.QA/9809156
\end{flushright}
\vskip.3in

\begin{center}
{\Large \bf Quasi-Hopf Superalgebras and Elliptic Quantum Supergroups}
\vskip.3in
{\large Yao-Zhong Zhang} and {\large Mark D. Gould}
\vskip.2in
{\em Department of Mathematics, University of Queensland, Brisbane,
     Qld 4072, Australia

Email: yzz@maths.uq.edu.au}
\end{center}

\vskip 2cm
\begin{center}
{\bf Abstract}
\end{center}
We introduce the quasi-Hopf superalgebras which are ${\bf Z}_2$ graded
versions of Drinfeld's quasi-Hopf algebras. We  describe the
realization of elliptic quantum supergroups as quasi-triangular
quasi-Hopf superalgebras obtained from twisting the normal quantum supergroups
by twistors which satisfy the graded shifted cocycle condition, 
thus generalizing the quasi-Hopf twisting procedure to the supersymmetric case. 
Two types of elliptic quantum supergroups are defined, that is the
face type ${\cal B}_{q,\lambda}({\cal G})$ and the vertex type
${\cal A}_{q,p}[\widehat{sl(n|n)}]$ (and 
${\cal A}_{q,p}[\widehat{gl(n|n)}]$), where ${\cal G}$ is any Kac-Moody
superalgebra with symmetrizable generalized Cartan matrix. It appears
that the vertex type twistor can be constructed only for
$U_q[\widehat{sl(n|n)}]$ in a non-standard system of simple roots, all
of which are fermionic.

\vskip 3cm
\noindent{\bf Mathematics Subject Classifications (1991):} 81R10, 17B37, 16W30

\end{titlepage}


\def\a{\alpha}
\def\b{\beta}
\def\d{\delta}
\def\e{\epsilon}
\def\ve{\varepsilon}
\def\g{\gamma}
\def\k{\kappa}
\def\l{\lambda}
\def\o{\omega}
\def\t{\theta}
\def\s{\sigma}
\def\D{\Delta}
\def\L{\Lambda}

\def\G{{\cal G}}
\def\hG{{\hat{\cal G}}}
\def\R{{\cal R}}
\def\hR{{\hat{\cal R}}}
\def\C{{\bf C}}
\def\P{{\bf P}}
\def\Z2{{{\bf Z}_2}}
\def\T{{\cal T}}
\def\H{{\cal H}}
\def\trho{{\tilde{\rho}}}
\def\tphi{{\tilde{\phi}}}
\def\tT{{\tilde{\cal T}}}
\def\uqsnh{{U_q[\widehat{sl(n|n)}]}}
\def\uqs1h{{U_q[\widehat{sl(1|1)}]}}


\def\beq{\begin{equation}}
\def\eeq{\end{equation}}
\def\bea{\begin{eqnarray}}
\def\eea{\end{eqnarray}}
\def\ba{\begin{array}}
\def\ea{\end{array}}
\def\no{\nonumber}
\def\lt{\left}
\def\rt{\right}
\newcommand{\bq}{\begin{quote}}
\newcommand{\eq}{\end{quote}}

\newtheorem{Theorem}{Theorem}
\newtheorem{Definition}{Definition}
\newtheorem{Proposition}{Proposition}
\newtheorem{Lemma}{Lemma}
\newtheorem{Corollary}{Corollary}
\newcommand{\proof}[1]{{\bf Proof. }
        #1\begin{flushright}$\Box$\end{flushright}}

\newcommand{\sect}[1]{\setcounter{equation}{0}\section{#1}}
\renewcommand{\theequation}{\thesection.\arabic{equation}}

\sect{Introduction\label{intro}}

One of the aims of this paper is to introduce $\Z2$ graded versions of 
Drinfeld's quasi-Hopf algebras \cite{Dri90}, which are referred to as
quasi-Hopf superalgebras. We then introduce elliptic quantum
supergroups, which are defined as quasi-triangular quasi-Hopf
superalgebras arising from twisting the normal quantum supergroups by
twistors which satisfy the graded shifted cocycle condition, 
thus generalizing Drinfeld's quasi-Hopf twisting
procedure \cite{Bab96,Fro97,Jim97,Arn97,Enr97}
to the supersymmetric case. 
We adopt the approach in \cite{Jim97}
and construct two types of twistors, i.e. the face type twistor
associated to any Kac-Moody superalgebra $\G$ with a
symmetrizable generalized Cartan matrix and the vertex type twistor 
associated to $\widehat{sl(n|n)}$ in a non-standard 
simple root system in which  all simple roots are odd (or fermionic).
It should be pointed out that the face type twistors for certain
classes of {\it non-affine} simple superalgebras were also constructed in
\cite{Arn97}. 

The elliptic quantum groups \cite{Fod94,Fel95}
are believed to provide the 
underlying algebraic structures for integrable models based on elliptic
solutions of the (dynamical) Yang-Baxter equation, such as Baxter's
8-vertex model \cite{Bax72}, the ABF model \cite{And84}
and their group theoretical generalizations \cite{Bel81}\cite{Jim88}.
The elliptic quantum supergroups described in this paper are
expected to play a similar role in supersymmetric integrable models
based on elliptic solutions \cite{Baz85,Deg91}
of the graded (dynamical) Yang-Baxter equation.

\sect{Quasi-Hopf Superalgebras}

\begin{Definition}\label{quasi-bi}: 
A $\Z2$ graded quasi-bialgebra is a $\Z2$ graded 
unital associative algebra $A$
over a field $K$  which is equipped with algebra homomorphisms $\e: 
A\rightarrow K$ (co-unit), $\D: A\rightarrow A\otimes A$ (co-product)
and an invertible homogeneous element $\Phi\in A\otimes A\otimes A$ 
(co-associator) satisfying 
\bea
&& (1\otimes\D)\D(a)=\Phi^{-1}(\D\otimes 1)\D(a)\Phi,~~
	\forall a\in A,\label{quasi-bi1}\\
&&(\D\otimes 1\otimes 1)\Phi \cdot (1\otimes 1\otimes\D)\Phi
	=(\Phi\otimes 1)\cdot(1\otimes\D\otimes 1)\Phi\cdot (1\otimes
	\Phi),\label{quasi-bi2}\\
&&(\e\otimes 1)\D=1=(1\otimes\e)\D,\label{quasi-bi3}\\
&&(1\otimes\e\otimes 1)\Phi=1.\label{quasi-bi4}
\eea
\end{Definition}
(\ref{quasi-bi2}), (\ref{quasi-bi3}) and (\ref{quasi-bi4}) imply
that $\Phi$ also obeys
\beq
(\e\otimes 1\otimes 1)\Phi=1=(1\otimes 1\otimes\e)\Phi.\label{e(phi)=1}
\eeq
The multiplication rule for the tensor products is $\Z2$ graded and is defined
for homogeneous elements $a,b,a',b'\in A$ by
\beq
(a\otimes b)(a'\otimes b')=(-1)^{[b][a']}\,(aa'\otimes bb'),
\eeq
where $[a]\in{\bf Z}_2$ 
denotes the grading of the element $a$.  

\begin{Definition}\label{quasi-hopf}: A quasi-Hopf superalgebra is
a $\Z2$ graded quasi-bialgebra $(A,\D,\e,\Phi)$ equipped with a
$\Z2$ graded algebra anti-homomorphism $S: A\rightarrow A$ (anti-pode) and 
canonical elements $\a,~\b\in A$ such that
\bea
&& m\cdot (1\otimes\a)(S\otimes 1)\D(a)=\e(a)\a,~~~\forall
    a\in A,\label{quasi-hopf1}\\
&& m\cdot (1\otimes\b)(1\otimes S)\D(a)=\e(a)\b,~~~\forall a\in A,
     \label{quasi-hopf2}\\
&& m\cdot (m\otimes 1)\cdot (1\otimes\b\otimes\a)(1\otimes S\otimes
     1)\Phi^{-1}=1,\label{quasi-hopf3}\\
&& m\cdot(m\otimes 1)\cdot (S\otimes 1\otimes 1)(1\otimes\a\otimes
     \b)(1\otimes 1\otimes S)\Phi=1.\label{quasi-hopf4}
\eea
\end{Definition}
Here $m$ denotes the usual product map on $A$: $m\cdot (a\otimes b)=ab,~
\forall a,b\in A$. Note that since $A$ is associative we have
$m\cdot(m\otimes 1)=m\cdot (1\otimes m)$.
For the homogeneous elements $a,b\in A$, the antipode satisfies
\beq
S(ab)=(-1)^{[a][b]}S(b)S(a),
\eeq
which extends to inhomogeneous elements through linearity.

Applying $\e$ to defintion (\ref{quasi-hopf3}, \ref{quasi-hopf4}) 
we obtain, in view
of (\ref{quasi-bi4}), $\e(\a)\e(\b)=1$. It follows that the canonical
elements $\a, \b$ are both even. By applying $\e$ to
(\ref{quasi-hopf1}), we have $\e(S(a))=\e(a),~\forall a\in A$.

In the following we show that the category of quasi-Hopf superalgebras
is invariant under a kind of gauge transformation. Let $(A,\D,\e,\Phi)$
be a qausi-Hopf superalgebra, with $\a,\b, S$ satisfying
(\ref{quasi-hopf1})-(\ref{quasi-hopf4}), and let $F\in A\otimes A$
be an invertible homogeneous element satisfying the co-unit properties
\beq
(\e\otimes 1)F=1=(1\otimes \e)F.\label{e(f)=1}
\eeq
It follows that $F$ is even. Throughout we set
\bea
&&\D_F(a)=F\D(a)F^{-1},~~~\forall a\in A,\label{twisted-d}\\
&&\Phi_F=(F\otimes 1)(\D\otimes
    1)F\cdot\Phi\cdot(1\otimes\D)F^{-1}(1\otimes F^{-1}).\label{twisted-phi}
\eea
\begin{Theorem}\label{t-quasi-hopf}:
$(A,\D_F,\e,\Phi_F)$ defined by (\ref{twisted-d}, 
\ref{twisted-phi}) together with
$\a_F,\b_F, S_F$ given by
\beq
S_F=S,~~~\a_F=m\cdot(1\otimes\a)(S\otimes 1)F^{-1},~~~
	 \b_F=m\cdot(1\otimes\b)(1\otimes S)F,\label{twisted-s-ab}
\eeq
is also a quasi-Hopf superalgebra. The element $F$ is referred to as
a twistor, throughout.
\end{Theorem}

The proof of this theorem is elementary. For demonstration we 
show in some details the proof of the anti-pode properties. 
Care has to be taken
of the gradings in tensor product multiplications and also in extending
the antipode to the whole algebra. First of all let us state 
\begin{Lemma}\label{L1}: For any elements $\eta\in A\otimes A$ and
$\xi\in A\otimes A\otimes A$,
\bea
&&m\cdot(1\otimes\a_F)(S\otimes 1)\eta=m\cdot(1\otimes \a)
       (S\otimes 1)(F^{-1}\eta),\label{L11}\\
&&m\cdot(1\otimes\b_F)(1\otimes S)\eta=m\cdot(1\otimes\b)
       (1\otimes S)(\eta F),\label{L12}\\
&&m\cdot(m\otimes 1)\cdot(1\otimes\b_F\otimes\a_F)
	(1\otimes S\otimes 1)\xi\no\\
&&~~~~=m\cdot(m\otimes 1)\cdot
	(1\otimes\b\otimes\a)(1\otimes S\otimes 1)
	[(1\otimes F^{-1})\cdot\xi\cdot(F\otimes 1)],\label{L13}\\
&&m\cdot(m\otimes 1)\cdot(S\otimes 1\otimes 1)(1\otimes\a_F\otimes
	\b_F)(1\otimes 1\otimes S)\xi\no\\
&&~~~~=m\cdot(m\otimes 1)\cdot
	(S\otimes 1\otimes 1)(1\otimes\a\otimes\b)(1\otimes 1\otimes
	S)\no\\
&&~~~~\cdot[(F^{-1}\otimes 1)\cdot\xi\cdot(1\otimes F)].\label{L14}
\eea
\end{Lemma}
\vskip.1in
\noindent {\em Proof}: Write
$F=f_i\otimes f^i$ and $F^{-1}=\bar{f}_i\otimes\bar{f}^i$. Here
and throughout, summation convention on repeated indices is
assumed. Then
(\ref{twisted-s-ab}) can be written as
\beq
\a_F=S(\bar{f}_i)\a\bar{f}^i,~~~~~~~\b_F=f_i\b S(f^i).\label{twisted-ab}
\eeq
Further write $\eta=\eta_k\otimes\eta^k$ and
  $\xi=\sum_ix_i\otimes y_i\otimes z_i$.
Then 
\bea
{\rm l.h.s.~of ~(\ref{L11})}&=&m\cdot (1\otimes S(\bar{f}_i)\a\bar{f}^i)
   (S(\eta_k)\otimes\eta^k)=m\cdot(S(\eta_k)\otimes S(\bar{f}_i)\a
   \bar{f}^i\eta^k)\no\\
   &=&S(\eta_k)S(\bar{f}_i)\a\bar{f}^i\eta^k
      =S(\bar{f}_i\eta_k)\a\bar{f}^i\eta^k
      \times (-1)^{[\eta_k][\bar{f}_i]},\no\\
{\rm r.h.s.~of~(\ref{L11})}&=&m\cdot(1\otimes\a)(S\otimes 1)
   (\bar{f}_i\eta_k\otimes
   \bar{f}^i\eta^k)\times(-1)^{[\bar{f}^i][\eta_k]}\no\\
  & =& S(\bar{f}_i\eta_k)\a\bar{f}^i\eta^k\times (-1)^{[\bar{f}_i][\eta_k]},\no
\eea
thus proving ({\ref{L11}). (\ref{L12}) can be proved
similarly. As for (\ref{L13}) we have:
\bea
{\rm l.h.s.~of~(\ref{L13})}&=&\sum_ix_i\b_FS(y_i)\a_Fz_i
   =\sum_ix_if_j\b S(f^j)S(y_i)S(\bar{f}_k)\a\bar{f}^kz_i\no\\
   &=&\sum_ix_if_j\b S(\bar{f}_ky_if^j)\a\bar{f}^kz_i
    \times(-1)^{[y_i]([f_j]+[\bar{f}_k])+[\bar{f}_k][f_j]},\no\\
{\rm r.h.s.~of~(\ref{L13})}&=&m\cdot(m\otimes 1)\cdot(1\otimes\b\otimes\a)
   (1\otimes S\otimes 1)\no\\
   & &\cdot\sum_i[x_if_j\otimes \bar{f}_ky_if^j\otimes\bar{f}^kz_i]
    \times(-1)^{[y_i]([f_j]+[\bar{f}_k])+[\bar{f}_k][f_j]},\no\\
   &=&\sum_ix_if_j\b S(\bar{f}_ky_if^j)\a\bar{f}^kz_i
    \times(-1)^{[y_i]([f_j]+[\bar{f}_k])+[\bar{f}_k][f_j]},\no
\eea
where we have used the fact that the element $F$ is even.
(\ref{L14}) is proved similarly.

Now let us prove the property (\ref{quasi-hopf1}) for $\a_F$ and
$\D_F$. We write, following Sweedler,
\beq
\D(a)=\sum_{(a)}a_{(1)}\otimes a_{(2)}.
\eeq
Then, in view of lemma \ref{L1},
\bea
m\cdot(1\otimes\a_F)(S\otimes 1)\D_F(a)&=&m\cdot(1\otimes\a)(S\otimes
   1)(F^{-1}\D_F(a)\no\\
&=&m\cdot(1\otimes\a)(S\otimes 1)(\D(a)F^{-1})\no\\
&=&m\cdot(1\otimes\a)\sum_{(a)}(S(a_{(1)}\bar{f}_i)\otimes a_{(2)}
   \bar{f}^i)\times (-1)^{[\bar{f}_i][a_{(2)}]}\no\\
&=&S(\bar{f}_i)\sum_{(a)}S(a_{(1)})\a a_{(2)}\bar{f}^i\times (-1)^{[\bar{f}_i]
  ([a_{(1)}]+[a_{(2)}])}\no\\
&=&S(\bar{f}_i)\sum_{(a)}S(a_{(1)})\a a_{(2)}\bar{f}^i
   \times(-1)^{[\bar{f}_i][a]}\no\\
&=&(-1)^{[\bar{f}_i][a]}S(\bar{f}_i)\;\sum_{(a)}S(a_{(1)})\a
   a_{(2)}\bar{f}^i\no\\
&\stackrel{(\ref{quasi-hopf1})}{=}&
   S(\bar{f}_i)\e(a)\a\bar{f}^i\times(-1)^{[\bar{f}_i][a]}\no\\
&=&S(\bar{f}_i)\e(a)\a\bar{f}^i\stackrel{(\ref{twisted-ab})}{=}\e(a)\a_F,
\eea
where we have used the fact that 
\beq
\e(a)=0,~~~~{\rm if~} [a]=1.\label{e(a)=0}
\eeq
The property (\ref{quasi-hopf2}) for
$\b_F$ and $\D_F$ is proved similarly. We then prove property 
(\ref{quasi-hopf3}), which reads in terms of the twisted objects
\beq
m\cdot (m\otimes 1)\cdot (1\otimes\b_F\otimes\a_F)(1\otimes S\otimes
     1)\Phi_F^{-1}=1.\label{extra1}
\eeq
Let us write
\beq
\Phi^{-1}=\sum_\nu \bar{X}_\nu\otimes\bar{Y}_\nu\otimes\bar{Z}_\nu.
\eeq
Then, in view of (\ref{L13}),
\bea
{\rm l.h.s.~of~(\ref{extra1})}&=&m\cdot(m\otimes 1)\cdot
   (1\otimes\b\otimes\a)(1\otimes S\otimes 1)[(1\otimes F^{-1})
   \Phi_F^{-1}(F\otimes 1)]\no\\
&=&m\cdot(m\otimes 1)\cdot
   (1\otimes\b\otimes\a)(1\otimes S\otimes 1)[(1\otimes\D)F\cdot
   \Phi^{-1}\cdot(\D\otimes 1)F^{-1}]\no\\
&=&m\cdot(m\otimes 1)\cdot
   (1\otimes\b\otimes\a)(1\otimes S\otimes 1)\no\\
& &\cdot   \sum_{\nu,(f),(\bar{f})}
   [f_i\bar{X}_\nu\bar{f}_{j(1)}\otimes f_{(1)}^i\bar{Y}_\nu
   \bar{f}_{j(2)}\otimes f^i_{(2)}\bar{Z}_\nu\bar{f}^j]\no\\
& &\times (-1)^{([\bar{X}_\nu]+[\bar{f}_{j(1)}])([f^i_{(1)}]+[f^i_{(2)}])
   +[\bar{Z}_\nu]([\bar{f}_{j(1)}]+[\bar{f}_{j(2)}])
   +[\bar{Y}_\nu]([\bar{f}_{j(1)}]+[f^i_{(2)}])
   +[f^i_{(2)}][\bar{f}_{j(2)}]}\no\\
&=&\sum_\nu f_i\bar{X}_\nu\sum_{(\bar{f})}\bar{f}_{j(1)}\b
   S(\bar{f}_{j(2)})S(\bar{Y}_\nu)\sum_{(f)}S(f_{(1)}^i)\a
   f^i_{(2)}\bar{Z}_\nu\bar{f}^j\no\\
& &\cdot (-1)^{([\bar{X}_\nu]+[\bar{Y}_\nu])([f^i_{(1)}]+[f^i_{(2)}])
   +([\bar{Y}_\nu]+[\bar{Z}_\nu])([\bar{f}_{j(1)}]+[\bar{f}_{j(2)}])
   +([f^i_{(1)}]+[f^i_{(2)}])([\bar{f}_{j(1)}]+[\bar{f}_{j(2)}])}\no\\
&=&\sum_\nu f_i\bar{X}_\nu\sum_{(\bar{f})}\bar{f}_{j(1)}\b
   S(\bar{f}_{j(2)})S(\bar{Y}_\nu)\sum_{(f)}S(f_{(1)}^i)\a
   f^i_{(2)}\bar{Z}_\nu\bar{f}^j\no\\
& &\cdot (-1)^{([\bar{X}_\nu]+[\bar{Y}_\nu])[f^i]
   +([\bar{Y}_\nu]+[\bar{Z}_\nu])[\bar{f}_j]
   +[f^i][\bar{f}_j]}\no\\
&=&\sum_\nu f_i\bar{X}_\nu\cdot (-1)^{([\bar{X}_\nu]+[\bar{Y}_\nu])[f^i]
   +([\bar{Y}_\nu]+[\bar{Z}_\nu])[\bar{f}_j]
   +[f^i][\bar{f}_j]}\no\\
& &\cdot \sum_{(\bar{f})}\bar{f}_{j(1)}\b
   S(\bar{f}_{j(2)})S(\bar{Y}_\nu)\sum_{(f)}S(f_{(1)}^i)\a
   f^i_{(2)}\bar{Z}_\nu\bar{f}^j\no\\
&\stackrel{(\ref{quasi-hopf1}, \ref{quasi-hopf2})}{=}&
   \sum_\nu f_i\bar{X}_\nu\e(\bar{f}_{j})\b
   S(\bar{Y}_\nu)\e(f^i)\a
   \bar{Z}_\nu\bar{f}^j
   \cdot(-1)^{([\bar{X}_\nu]+[\bar{Y}_\nu])[f^i]
   +([\bar{Y}_\nu]+[\bar{Z}_\nu])[\bar{f}_j]
   +[f^i][\bar{f}_j]}\no\\
&\stackrel{(\ref{e(a)=0})}{=}&
   \sum_\nu f_i\bar{X}_\nu\e(\bar{f}_{j})\b
   S(\bar{Y}_\nu)\e(f^i)\a
   \bar{Z}_\nu\bar{f}^j\no\\
&=&m\cdot(m\otimes 1)\cdot
   (1\otimes\b\otimes\a)(1\otimes S\otimes 1)\no\\
& &\cdot   [((1\otimes\e)F\otimes 1)\cdot
   \Phi^{-1}\cdot((\e\otimes 1)F^{-1}\otimes 1)]\no\\
&\stackrel{(\ref{e(f)=1})}{=}&m\cdot(m\otimes 1)\cdot
   (1\otimes\b\otimes\a)(1\otimes S\otimes 1)\Phi^{-1}
   \stackrel{(\ref{quasi-hopf3})}{=}1.\no
\eea
The property (\ref{quasi-hopf4}) for the twisted objects, which reads,
\beq
m\cdot(m\otimes 1)\cdot (S\otimes 1\otimes 1)(1\otimes\a_F\otimes
     \b_F)(1\otimes 1\otimes S)\Phi_F=1,
\eeq
is proved in a similar way.

\begin{Definition}\label{quasi-quasi}: A quasi-Hopf
superalgebra $(A,\D,\e,\Phi)$ is called quasi-triangular if there
exists an invertible homogeneous element $\R\in A\otimes A$ such that
\bea
&&\D^T(a)\R=\R\D(a),~~~~\forall a\in A,\label{dr=rd}\\
&&(\D\otimes 1)\R=\Phi^{-1}_{231}\R_{13}\Phi_{132}\R_{23}\Phi^{-1}_{123},
   \label{d1r}\\
&&(1\otimes \D)\R=\Phi_{312}\R_{13}\Phi^{-1}_{213}\R_{12}\Phi_{123}.
   \label{1dr}
\eea
\end{Definition}
Throughout, $\D^T=T\cdot\D$ with $T$ being the graded twist map
which is defined, for homogeneous elements $a,b\in A$, by
\beq
T(a\otimes b)=(-1)^{[a][b]}b\otimes a;
\eeq
and $\Phi_{132}$ {\it etc} are derived from $\Phi\equiv\Phi_{123}$
with the help of $T$
\bea
&&\Phi_{132}=(1\otimes T)\Phi_{123},\no\\
&&\Phi_{312}=(T\otimes 1)\Phi_{132}=(T\otimes 1)
   (1\otimes T)\Phi_{123},\no\\
&&\Phi^{-1}_{231}=(1\otimes T)\Phi^{-1}_{213}=(1\otimes T)
   (T\otimes 1)\Phi^{-1}_{123},\no
\eea
and so on. We remark that our convention differs from the usual one
which employs the inverse permutation on the positions (c.f. 
\cite{Jim97}).

It is easily shown that the properties (\ref{dr=rd})-(\ref{1dr})
imply the graded Yang-Baxter type equation,
\beq
\R_{12}\Phi^{-1}_{231}\R_{13}\Phi_{132}\R_{23}\Phi^{-1}_{123}
  =\Phi^{-1}_{321}\R_{23}\Phi_{312}\R_{13}\Phi^{-1}_{213}\R_{12},
  \label{quasi-ybe}
\eeq
which is referred to as the graded  quasi-Yang-Baxter equation, 
and the co-unit properties of $\R$:
\beq
(\e\otimes 1)\R=1=(1\otimes \e)\R.\label{e(R)=1}
\eeq

\begin{Theorem}\label{t-quasi-quasi}: Denoting by the set
$(A,\D,\e,\Phi,\R)$  a
quasi-triangular quasi-Hopf superalgebra, then $(A, \D_F, \e, \Phi_F, \R_F)$
is also a quasi-triangular quasi-Hopf superalgebra, with the choice of
$R_F$ given by
\beq
\R_F=F^T \R F^{-1},\label{twisted-R}
\eeq
where $F^T=T\cdot F\equiv F_{21}$. Here $\D_F$ and $\Phi_F$ are given
by (\ref{twisted-d}) and (\ref{twisted-phi}), respectively.
\end{Theorem}

The proof of this theorem is elementary
computation. As an example,
let us illustrate the proof of the property (\ref{d1r}) for $\D_F,
\R_F$ and $\Phi_F$. Applying the homomorphism $T\otimes 1$ to
$(\Phi^{-1}_F)_{123}$, one obtains
\bea
(\Phi_F^{-1})_{213}&=&F_{13}(T\otimes 1)(1\otimes \D)F\cdot
   \Phi^{-1}_{213}\cdot(\D^T\otimes)F^{-1}\cdot (F^T)^{-1}_{12}\no\\
&=&F_{13}\sum_{(f)}(-1)^{[f^i_{(1)}][f_i]}(f^i_{(1)}\otimes f_i\otimes
   f^i_{(2)})\Phi^{-1}_{213}(\D^T\otimes 1)F^{-1}\cdot (F^T)^{-1}_{12},
   \label{1}
\eea
which gives rise to, by applying
the homomorphism $1\otimes T$ to both sides,
\bea
(\Phi^{-1}_F)_{231}&=&F_{12}\sum_{(f)}(-1)^{([f^i_{(1)}]+[f^i_{(2)}])[f_i]}
   (f^i_{(1)}\otimes 
   f^i_{(2)}\otimes f_i)\Phi^{-1}_{231}(1\otimes T)(\D^T\otimes 1)F^{-1}
   \cdot (F^T)^{-1}_{13}\no\\
&=&F_{12}(\D\otimes 1)F^T\cdot\Phi^{-1}_{231}(1\otimes T)(\D^T\otimes
   1)F^{-1}\cdot (F^T)^{-1}_{13}.\label{2}
\eea
Then,
\bea
(\D_F\otimes 1)\R_F&=&(F\otimes 1)(\D\otimes 1)\R_F\cdot(F^{-1}\otimes
   1)\no\\
&=&F_{12}(\D\otimes 1)(F^T\R F^{-1})\cdot F^{-1}_{12}\no\\
&=&F_{12}(\D\otimes 1)F^T(\D\otimes 1)\R(\D\otimes 1)F^{-1}\cdot
   F^{-1}_{12}\no\\
&\stackrel{(\ref{d1r})}{=}&F_{12}(\D\otimes 1)F^T\cdot\Phi^{-1}_{231}\R_{13}
  \Phi_{132}\R_{23}\Phi^{-1}_{123}(\D\otimes 1)F^{-1}\cdot
  F^{-1}_{12}\no\\
&\stackrel{(\ref{2})}{=}&(\Phi^{-1}_F)_{231}(F^T)_{13}(1\otimes T)
  (\D^T\otimes 1)F\cdot \R_{13}
  \Phi_{132}\R_{23}\Phi^{-1}_{123}(\D\otimes 1)F^{-1}\cdot
  F^{-1}_{12}\no\\
&\stackrel{(\ref{twisted-phi})}{=}&(\Phi^{-1}_F)_{231}(F^T)_{13}(1\otimes T)
  (\D^T\otimes 1)F\no\\
& &  \cdot \R_{13}
  \Phi_{132}\R_{23}(1\otimes \D)F^{-1}\cdot
  F^{-1}_{23}(\Phi^{-1}_F)_{123}\no\\
&=&(\Phi^{-1}_F)_{231}(F^T)_{13}(1\otimes T)[
  (\D^T\otimes 1)F\cdot \R_{12}]\no\\
& &\cdot  \Phi_{132}\R_{23}(1\otimes \D)F^{-1}\cdot 
  F^{-1}_{23}(\Phi^{-1}_F)_{123}\no\\
&\stackrel{(\ref{dr=rd})}{=}&(\Phi^{-1}_F)_{231}(F^T)_{13}(1\otimes T)[
  \R_{12}(\D\otimes 1)F]\no\\
& &\cdot  \Phi_{132}(1\otimes \D^T)F^{-1}\cdot \R_{23}
  F^{-1}_{23}(\Phi^{-1}_F)_{123}\no\\
&=&(\Phi^{-1}_F)_{231}(F^T)_{13}\R_{13}(1\otimes T)[
  (\D\otimes 1)F]\no\\
& &\cdot  \Phi_{132}(1\otimes \D^T)F^{-1}\cdot \R_{23}
  F^{-1}_{23}(\Phi^{-1}_F)_{123}\no\\
&\stackrel{(\ref{twisted-R})}{=}&(\Phi^{-1}_F)_{231}(\R_F)_{13}F^{-1}_{13}
  (1\otimes T)[(\D\otimes 1)F]\no\\
& &\cdot  \Phi_{132}(1\otimes \D^T)F^{-1}
  (F^T)^{-1}_{23}(\R_F)_{23}(\Phi^{-1}_F)_{123}\no\\
&=&(\Phi^{-1}_F)_{231}(\R_F)_{13}(1\otimes T)[F^{-1}_{12}
  (\D\otimes 1)F
  \Phi_{123}(1\otimes \D)F^{-1}\cdot
  F^{-1}_{23}]\no\\
& &\cdot  (\R_F)_{23}(\Phi^{-1}_F)_{123}\no\\
&\stackrel{(\ref{twisted-phi})}{=}&(\Phi^{-1}_F)_{231}(\R_F)_{13}(1\otimes T)
  (\Phi_F)_{123}\cdot
  (\R_F)_{23}(\Phi^{-1}_F)_{123}\no\\
&=&(\Phi^{-1}_F)_{231}(\R_F)_{13}(\Phi_F)_{132}
  (\R_F)_{23}(\Phi^{-1}_F)_{123}.
\eea

\vskip.1in
Let us now consider the special case that $A$ arises from a normal
quasi-triangular Hopf superalgebra via twisting with $F$. 
A quasi-triangular Hopf superalgebra is a quasi-triangular 
quasi-Hopf superalgebra with $\a=\b=1,~\Phi=1\otimes 1\otimes 1$. 
Hence $A$ has the following $\Z2$ graded quasi-Hopf algebra structure,
\bea
&&\D_F(a)=F\D(a)F^{-1},~~~~\forall a\in A,\no\\
&&\Phi_F=F_{12}\cdot(\D\otimes 1)F\cdot(1\otimes\D)F^{-1}\cdot 
   F_{23}^{-1},\no\\
&&\a_F=m\cdot(S\otimes 1)F^{-1},~~~~\b_F=m\cdot(1\otimes S)F,\no\\
&&\R_F=F^T \R F^{-1}.\label{twisting-qg}
\eea
The twisting procedure is particularly interesting
when the twistor $F\in A\otimes A$ depends on an element $\l\in A$,
i.e. $F=F(\l)$, and is a shifted cocycle in the following sense.
Here $\l$ is assumed to depend on one (or possible several) parameters.

\begin{Definition}: A twistor $F(\l)$ depending
on $\l\in A$ is a shifted cocycle if it satisfies the graded
shifted cocycle condition:
\beq
F_{12}(\l)\cdot (\D\otimes 1)F(\l)=F_{23}(\l+h^{(1)})\cdot
   (1\otimes \D)F(\l),\label{shifted-cocycle}
\eeq
where $h^{(1)}=h\otimes 1\otimes 1$ and $h\in A$ is fixed.
\end{Definition}

Let $(A,\D_\l,\e,\Phi(\l),\R(\l))$ be the quasi-triangular quasi-Hopf
superalgebra obtained from twisting the quasi-triangular Hopf superalgebra
by the twistor $F(\l)$. Then

\begin{Proposition}: We have
\bea
&&\Phi(\l)\equiv\Phi_F=F_{23}(\l+h^{(1)})F_{23}(\l)^{-1},\label{phi-lambda}\\
&&\D_\l(a)^T \R(\l)=\R(\l)\D_\l(a),~~~~\forall a\in
    A,\label{d-lambda}\\
&&(\D_\l\otimes 1)\R(\l)=\Phi_{231}(\l)^{-1}\R_{13}(\l)\R_{23}
   (\l+h^{(1)}),\label{d-lambda-1-r}\\
&&(1\otimes \D_\l)\R(\l)=\R_{13}(\l+h^{(2)})\R_{12}(\l)\Phi_{123}(\l).
   \label{1-d-lambda-r}
\eea
As a corollary, $\R(\l)$ satisfies the graded dynamical Yang-Baxter
equation
\beq
\R_{12}(\l+h^{(3)})\R_{13}(\l)\R_{23}(\l+h^{(1)})
  =\R_{23}(\l)\R_{13}(\l+h^{(2)})\R_{12}(\l).\label{dynamical-ybe}
\eeq
\end{Proposition}

\sect{Quantum Supergroups}

Let $\G$ be a Kac-Moody superalgebra \cite{Kac77,Kac78}
with a symmetrizable generalized
Cartan matrix $A=(a_{ij})_{i,j,\in I}$. As is well-known, a given Kac-Moody
superalgebra allows many inequivalent systems of simple roots.
A system of simple roots is called distinguished if it has
minimal odd roots.
Let $\{\a_i,~i\in I\}$ denote a
chosen set of simple roots. Let $(~,~)$ be a fixed 
invariant bilinear form on the root space of $\G$. Let $\H$ be
the Cartan subalgebra and throughout we identify the
dual $\H^*$ with $\H$ via $(~,~)$.
The generalized Cartan matrix $A=(a_{ij})_{i,j\in I}$ is
defined from the simple roots by
\beq
a_{ij}=\left \{
\begin{array}{l}
\frac{2(\a_i,\a_j)}{(\a_i,\a_i)},~~~~{\rm if}~(\a_i,\a_i)\neq 0\\
(\a_i,\a_j),~~~~{\rm if}~(\a_i,\a_i)=0
\end{array}
\right .
\eeq

As we mentioned in the previous section,
quantum Kac-Moody superalgebras are quasi-triangular quasi-Hopf
superalgebras with $\a=\b=1,~\Phi=1\otimes 1\otimes 1$.
We shall not give the standard relations obeyed by the simple 
generators (or Chevalley generators) $\{h_i,~e_i,~f_i,~i\in I\}$ of
$U_q({\cal G})$, but mention that for certain types of Dynkin diagrams
extra $q$-Serre relations are
needed in the defining relations. We adopt the following
graded Hopf algebra structure
\bea
\D(h)&=&h\otimes 1+1\otimes h,\no\\
\D(e_i)&=&e_i\otimes 1+t_i\otimes e_i,~~~~
\D(f_i)=f_i\otimes t_i^{-1}+1\otimes f_i,\no\\
\e(e_i)&=&\e(f_i)=\e(h)=0,\no\\
S(e_i)&=&-t_i^{-1}e_i,~~~~S(f_i)=-f_it_i,~~~~S(h)=-h,\label{e-s} 
\eea
where $i\in I$, $t_i=q^{h_i}$ and $h\in {\cal H}$.

The canonical element $\R$ is called the universal R-matrix of $U_q(\G)$, which 
satisfies the basic properties (e.g. (\ref{dr=rd})-(\ref{1dr}) with
$\Phi=1\otimes 1\otimes 1$ and (\ref{e(R)=1}}))
\bea
&&\D^T(a)\R=\R\D(a),~~~~\forall a\in U_q(\G),\no\\
&&(\D\otimes 1)\R=\R_{13}\R_{23},\no\\
&&(1\otimes\D)\R=\R_{13}\R_{12},\no\\
&&(\e\otimes 1)\R=(1\otimes\e)\R=1,\label{D-hR}
\eea
and the graded Yang-Baxter equation (c.f. (\ref{quasi-ybe}) with
$\Phi=1\otimes 1\otimes 1$)
\beq
\R_{12}\R_{13}\R_{23}=\R_{23}\R_{13}\R_{12}.
\eeq

The Hopf superalgebra $U_q(\G)$ contains two important Hopf subalgebras
$U_q^+$ and $U_q^-$ which are generated by
$e_i$ and
$f_i$, respectively. By Drinfeld's 
quantum double construction, the
universal R-matrix  $\R$ can be written in the form 
\beq
\R=\lt(1\otimes 1+\sum_t\,a^t\otimes a_t\rt)\cdot q^{-\T}
\eeq
where $\{a^t\}\in U_q^+,~\{a_t\}\in U_q^-$. The element $\T$
is defined as follows.
If the symmetrical Cartan matrix is non-degenerate,
then $\T$ is the usual canonical element of $\H\otimes\H$.
Let $\{h_l\}$ be a basis of $\H$ and
$\{h^l\}$ be its dual basis. Then $\T$ can be written as
\beq
\T=\sum_lh_l\otimes h^l.\label{T}
\eeq
In the case of a degenerate symmetrical
Cartan matrix, we extend the Cartan subalgebra $\H$ by adding some elements
to it in such a way that the extended symmetrical Cartan matrix is
non-degenerate \cite{Kho91}. Then $\T$ stands for the canonical element
of the extended Cartan subalgebra. It still takes the
form (\ref{T}) but now $\{h_l\}~(\{h^l\})$ is understood to be the
(dual) basis of the extended Cartan subalgebra. After such enlargement,
one has $h=\sum_l(h^l,h)h_l=\sum_l(h_l,h)h^l$ for any given $h$
in the enlarged Cartan subalgebra. 

For later use, we work out the explicit form of the
universal R-matrix for the
simplest quantum affine superalgebra $\uqs1h$. This algebra is generated
by Chevalley generators $\{e_i, f_i, h_i, d, i=0,1\}$
with $e_i,~f_i$ odd and $h_i,~d$ even. Here and throughout $d$ stands
for the derivation operator. Let us write $h_i=\a_i$.
Then we have $h_0=\d-\ve_1+\d_1,~h_1=\ve_1-\d_1$, where
$\{\ve_1,\d_1,\d\}$ satisfy 
$(\ve_1,\ve_1)=1
=-(\d_1,\d_1),~(\ve_1,\d_1)=(\d,\d)=(\d,\ve_1)=(\d,\d_1)=0$.
We extend the Cartan
subalgebra by adding to it the element $h_{\rm ex}=\ve_1+\d_1$.
A basis for the enlarged Cartan subalgebra is thus
$\{h_{\rm ex},h_0,h_1,d\}$.
It is easily shown that the dual basis
is $\{h^{\rm ex},h^0,h^1,c\}$, where $h^{\rm ex}=\frac{1}{2}(\ve_1-\d_1)=
\frac{1}{2}h_1,~
h^0=d,~h^1=\ve_1+d-\frac{1}{2}(\ve_1-\d_1)=d+\frac{1}{2}h_{\rm ex}$.
As is well-known, $\uqs1h$ can also be realized in terms of the Drinfeld
generators \cite{Dri88} 
$\{X^\pm_n, H_n, H^{\rm ex}_n, n\in {\bf Z},c,d\}$, where
$X^\pm_n$ are odd and all other generators are even. The relations 
satisfied by the Drinfeld generators read \cite{Zha97}
\bea
&&[c,a]=[H_0,a]=[d,d]=[H_n,H_m]=[H^{\rm ex}_n,H^{\rm ex}_m]=0,~~~
   \forall a\in \uqs1h,\no\\
&&q^{H^{\rm ex}_0}X^\pm_nq^{-H^{\rm ex}_0}=q^{\pm 2}X^\pm_n,\no\\
&&[d,X^\pm_n]=nX^\pm_n,~~~[d,H_n]=nH_n,~~~[d,H^{\rm ex}_n]=nH^{\rm
   ex}_n,\no\\
&&[H_n, H_m^{\rm ex}]=\d_{n+m,0}\frac{[2n]_q[nc]_q}{n},\no\\
&&[H^{\rm ex}_n,
   X^\pm_m]=\pm\frac{[2n]_q}{n}X^\pm_{n+m}q^{\mp|n|c/2},\no\\
&&[H_n,X^\pm_m]=0=[X^\pm_n,X^\pm_m],\no\\
&&[X^+_n, X^-_m]=\frac{1}{q-q^{-1}}\lt(q^{\frac{c}{2}(n-m)}
  \psi^+_{n+m}-q^{-\frac{c}{2}(n-m)}\psi^-_{n+m}\rt),\label{Drinfeld-sl11}
\eea
where $[x]_q=(q^x-q^{-x})/(q-q^{-1})$,
$[a,b]\equiv ab-(-1)^{[a][b]}ba$ denotes the supercommutator
and $\psi^\pm_{\pm n}$ are related to $H_{\pm n}$ by relations
\beq
\sum_{n\geq 0}\psi^\pm_{\pm n}z^{\mp n}=q^{\pm H_0}\exp\lt(
  \pm(q-q^{-1})\sum_{n>0}H_{\pm n}z^{\mp n}\rt).
\eeq
The relationship between the Drinfeld generators and the 
Chevalley generators is
\bea
&&e_1=X^+_0,~~~~f_1=X^-_0,~~~~h_1=H_0,~~~~h_{\rm ex}=H^{\rm ex}_0,\no\\
&&e_0=X^-_1 q^{-H_0},~~~~f_0=-q^{H_0}X^+_{-1},~~~~h_0=c-H_0.
\eea
With the help of the Drinfeld generators, we find
the following universal R-matrix 
\beq
\R=\R'\cdot q^{-\T},
\eeq
where
\bea
\T&=&h_{\rm ex}\otimes h^{\rm ex}+h_0\otimes h^0+h_1\otimes h^1
  +d\otimes c\no\\
  &=&\frac{1}{2}(H_0\otimes H_0^{\rm ex}
  +H^{\rm ex}_0\otimes H_0)+c\otimes d+d\otimes c,\no\\
\R'&=&\R^< \,\R^0 \,\R^>,\no\\
\R^<&=&\prod^{\rightarrow}_{n\geq 0}\exp\lt[(q-q^{-1})
  (q^{-nc/2}X^+_n\otimes q^{nc/2}X^-_{-n})\rt],\no\\
\R^0&=&\exp\lt[-(q-q^{-1})\sum_{n=1}^\infty\frac{n}{[2n]_q}
  (H_n\otimes H^{\rm ex}_{-n}+H^{\rm ex}_n\otimes H_{-n})\rt],\no\\
\R^>&=&\prod^{\leftarrow}_{n\geq 0}\exp\lt[-(q-q^{-1})
  (X^-_{n+1}q^{nc/2-H_0}\otimes q^{-nc/2+H_0}X^+_{-n-1})\rt].
  \label{sl11-R}
\eea
Here and throughout, 
\beq
\prod_{k\geq 0}^{\rightarrow}A_k=A_0A_1A_2\cdots,~~~~~
\prod_{k\geq 0}^{\leftarrow}A_k=\cdots A_2A_1A_0.
\eeq
It seems to us that even for this simplest quantum affine superalgebra
$\uqs1h$ the universal R-matrix has not been written down
in its explicit form before.

Let us compute the image of $\R$ in the 2-dimensional evaluation
representaion $(\pi,V)$ of $\uqs1h$, where $V=\C^{1|1}=\C v_1\oplus\C
v_2$ with $v_1$ even and $v_2$ odd.
Let $e_{ij}$ be the $2\times 2$ matrix whose $(i,j)$-element is unity
and zero otherwise. 
In the homogeneous gradation, the simple 
generators are represented by
\bea
&&e_1=\sqrt{[\t]_q}e_{12},~~~f_1=\sqrt{[\t]_q}e_{21},~~~
  h_1=\t(e_{11}+e_{22}),~~~h_{\rm ex}=2e_{11}+c_0(e_{11}+e_{22}),\no\\
&&e_0=z\sqrt{[\t]_q}e_{21},~~~f_0=-z^{-1}\sqrt{[\t]_q}e_{12},~~~
  h_0=-\t(e_{11}+e_{22}),\label{rep-V}
\eea
where $\t$ and $c_0$ are arbitrary constants. Then it can be shown
that the Drinfeld generators are represented by
\bea
&&H_n=z^n\frac{[n\t]_q}{n}(e_{11}+e_{22}),~~~~
  H^{\rm ex}_n=z^n\frac{[2n]_q}{n}q^{n\t}e_{11}+z^nc_n(e_{11}+e_{22}),\no\\
&&X^+_n=z^nq^{n\t}\sqrt{[\t]_q}e_{12},~~~~
  X^-_n=z^nq^{n\t}\sqrt{[\t]_q}e_{21},
\eea
where again $c_n$ are arbitrary constants. In the following we set
$c_n$ to be zero. Then the image $R_{VV}(z;\t,\t')=(\pi_\t\otimes\pi_{\t'})\R$ 
depends on two extra non-additive parameters $\t,~\t'$, and is given by
\bea
R_{VV}(z;\t,\t')&=& \frac{q^{-\t-\t'}-z}
   {1-zq^{-\t-\t'}}e_{11}\otimes e_{11}+e_{22}\otimes e_{22}
   +\frac{q^{-\t'}-zq^{-\t}}
   {1-zq^{-\t-\t'}}e_{11}\otimes e_{22}\no\\
& &   +\frac{q^{-\t}-zq^{-\t'}}
   {1-zq^{-\t-\t'}}e_{22}\otimes e_{11}
   +\sqrt{[\t]_q[\t']_q}q^{-\t}\frac{q-q^{-1}}
   {1-zq^{-\t-\t'}}e_{12}\otimes e_{21}\no\\
& &   -\sqrt{[\t]_q[\t']_q}q^{-\t'}\frac{z(q-q^{-1})}
   {1-zq^{-\t-\t'}}e_{21}\otimes e_{12}. \label{sl11-r}
\eea
(\ref{sl11-r}) is nothing but the R-matrix obtained in \cite{Bra94} by
solving the Jimbo equation.

\sect{Elliptic Quantum Supergroups}

Following Jimbo et al \cite{Jim97}, we define elliptic quantum supergroups
to be quasi-triangular quasi-Hopf superalgebras obtained from
twisting the normal quantum supergroups (which are quasi-triangular quasi-Hopf
superalgebras with $\a=\b=1,~\Phi=1\otimes 1\otimes 1$)
by twistors which satisfy the graded 
shifted cocycle condition.

\subsection{Elliptic Quantum Supergroups of Face Type}

Let $\rho$ be an element in the (extended) Cartan subalgebra such
that $(\rho,\a_i)=(\a_i,\a_i)/2$ for all $i\in I$, and
\beq
\phi={\rm Ad}(q^{\frac{1}{2}\sum_lh_lh^l-\rho}),
\eeq
be an automorphism of $U_q(\G)$. Here
$\{h_l\},~\{h^l\}$ are as in (\ref{T}) and are
the dual basis of the (extended) Cartan
subalgebra. Namely,
\beq
\phi(e_i)=e_it_i,~~~~\phi(f_i)=t_i^{-1}f_i,~~~~\phi(q^h)=q^h.
\eeq
In the following we consider the special case in which
the element $\l$ introduced before
belongs to the (extended) Cartan subalgebra. Let
\beq
\phi_\l=\phi^2\cdot{\rm Ad}(q^{2\l})={\rm Ad}(q^{\sum_lh_lh^l-2\rho+
   2\l})
\eeq
be an automorphism depending on the element $\l$ and $\R$ be the
universal R-matrix of $U_q(\G)$. Following Jimbo et al \cite{Jim97},
we define a twistor $F(\l)$ by the infinite product
\beq
F(\l)=\prod_{k\geq 1}^{\leftarrow}\lt(\phi_\l^k\otimes 1\rt)
  \lt(q^\T\R\rt)^{-1}.\label{twistor-f}
\eeq
It is easily seen that $F(\l)$ is a formal series in parameter(s)
in $\l$ with leading term 1. Therefore the infinite product makes
sense. The twistor $F(\l)$ is referred to as face type twistor.
It can be shown that $F(\l)$ satisfies the graded shifted cocycle condition
\beq
F_{12}(\l)(\D\otimes 1)F(\l)=F_{23}(\l+h^{(1)})(1\otimes\D)F(\l),
  \label{shifted-f}
\eeq
where, if $\l=\sum_l\l_lh^l$, then $\l+h^{(1)}=\sum_l(\l_l+h_l^{(1)})
h^l$. 
The proof of (\ref{shifted-f}) is identical to the non-super case given by
Jimbo et al \cite{Jim97}, apart from the use of the graded tensor products. 
Moreover, it is easily seen that $F(\l)$ obeys the co-unit property
\beq
(\e\otimes 1)F(\l)=(1\otimes\e)F(\l)=1.
\eeq
We have

\begin{Definition} {\bf (Face type elliptic quantum supergroup)}: We define
 elliptic quantum supergroup ${\cal B}_{q,\l}(\G)$ of face type to be
the quasi-triangular quasi-Hopf superalgebra
$(U_q(\G),\D_\l,\e,\Phi(\l),\R(\l))$ together with the graded algebra
anti-homomorphism $S$ defined by (\ref{e-s}) 
and $\a_\l=m\cdot(S\otimes 1)F(\l)^{-1}$,
$\b_\l=m\cdot(1\otimes S)F(\l)$. Here $\e$ is defined by (\ref{e-s}),
and
\bea
&&\D_\l(a)=F(\l)\D(a)F(\l)^{-1},~~~\forall a\in U_q(\G),\no\\
&&\R(\l)=F(\l)^T\R F(\l)^{-1},\no\\
&&\Phi(\l)=F_{23}(\l+h^{(1)})F_{23}(\l)^{-1}.
\eea
\end{Definition}

We now consider the particularly interesting case where $\G$ is of affine type.
Then $\rho$ contains two parts
\beq
\rho=\bar{\rho}+gd,
\eeq
where $g=(\psi,\psi+2\bar{\rho})/2$, $\bar{\rho}$ is the graded half-sum
of positive roots of the non-affine part $\bar{\G}$ and $\psi$ is  highest
root of $\bar{\G}$; $d$ is the derivation operator which gives the
homogeneous gradation
\beq
[d, e_i]=\d_{i0}e_i,~~~~[d,f_i]=-\d_{i0}f_i,~~~~i\in I.
\eeq
We also set
\beq
\l=\bar{\l}+(r+g)d+s'c,~~~~r,s'\in \C,
\eeq
where $\bar{\l}$ stands for the projection of $\l$ onto the
(extended) Cartan subalgebra of $\bar{\G}$.
Denoting by $\{\bar{h}_j\},~\{\bar{h}^j\}$ the dual basis of the
(extended) Cartan subalgebra of $\bar{\G}$
and setting $p=q^{2r}$, we can decompose $\phi_\l$ into two parts
\beq
\phi_\l={\rm Ad}(p^dq^{2cd})\cdot\bar{\phi}_\l,~~~~
   \bar{\phi}_\l={\rm
   Ad}(q^{\sum_j\bar{h}_j\bar{h}^j+2(\bar{\l}-\bar{\rho})}).
\eeq
Introduce a formal parameter $z$ (which will
be identified with spectral parameter) into $\R$ and $F(\l)$ by setting
\bea
&&\R(z)={\rm Ad}(z^d\otimes 1)\R,\no\\
&&F(z,\l)={\rm Ad}(z^d\otimes 1)F(\l),\no\\
&&\R(z,\l)={\rm Ad}(z^d\otimes 1)\R(\l)=F(z^{-1},\l)^T
   \R(z)F(z,\l)^{-1}.
\eea
Then it can be shown from the definition of $F(\l)$
that $F(z,\l)$ satisfies the difference equation
\bea
&&F(pq^{2c^{(1)}}z,\l)=(\bar{\phi}_\l\otimes 1)^{-1}(F(z,\l))\cdot
  q^\T \R(pq^{2c^{(1)}}z),\no\\
&&F(0,\l)=F_{\bar{\G}}(\bar{\l}).\label{dif-face}
\eea
The initial condition follows from the fact that $\R(z)q^{d\otimes c
+c\otimes d}|_{z=0}$ reduces to the universal R-matrix of $U_q(\bar{\G})$.

Let us give some examples.
\vskip.1in
\noindent\underline{\bf The case ${\cal B}_{q,\l}[sl(1|1)]$}:
\vskip.1in
In this case the universal R-matrix is given simply by
\bea
&&\R=\exp[(q-q^{-1})\,e\otimes f]\;q^{-\T}=[1+(q-q^{-1})\,e\otimes f]
   \;q^{-\T},\no\\
&&\T=\frac{1}{2}(h\otimes h_{\rm ex}+h_{\rm ex}\otimes h).
\eea
Let us write
\beq
\l=(s'+1)\frac{1}{2}h+s\frac{1}{2}h_{\rm ex},~~~~s',s\in \C.
\eeq
Since $h$ commutes with everything, $\phi_\l$ is independent of $s'$.
Set $w=q^{2(s+h)}$, we have
\beq
\phi_\l={\rm Ad}(w^{\frac{1}{2}h_{\rm ex}}).
\eeq
The formula for the twistor becomes 
\bea
F(w)&=&\prod_{k\geq 1}\lt(1-(q-q^{-1})w^k\,q^{-h}e\otimes fq^h\rt)\no\\
    &=&1-(q-q^{-1})\sum_{k=1}^\infty w^k\,q^{-h}e\otimes fq^h\no\\
    &=&1-(q-q^{-1})\frac{w}{1-w}\,q^{-h}e\otimes fq^h.
\eea

\vskip.1in
\noindent\underline{\bf The case ${\cal B}_{q,\l}[\widehat{sl(1|1)}]$}:
\vskip.1in

Taking a basis $\{c,d,h,h_{\rm ex}\}$ of the enlarged Cartan subalgebra
of $\widehat{sl(1|1)}$, we write
\beq
\l=rd+s'c+(s''+1)\frac{1}{2}h+s\frac{1}{2}h_{\rm ex},~~~~r,s',s'',s\in \C.
\eeq
Then $\phi_\l$ is independent of $s'$ and $s''$. Set
\beq
p=q^{2r},~~~~~~w=q^{2(s+h)},
\eeq
Set $F(z;p,w)\equiv F(z,\l)$. Then (\ref{dif-face}) take the form
\bea
&&F(pq^{2c^{(1)}}z;p,w)=(\bar{\phi}_w^{-1}\otimes 1)(F(z;p,w))\cdot
  q^\T \R(pq^{2c^{(1)}}z),\label{dif-face-sl11}\\
&&F(0;p,w)=F_{sl(1|1)}(w),\label{initial-sl11-face}
\eea
where $\bar{\phi}_w={\rm Ad}(w^{\frac{1}{2}h_{\rm ex}})$. 

The image of (\ref{dif-face-sl11}) in the two-dimensional representation 
$(\pi, V)$ given by (\ref{rep-V}) (by setting $\t=1$)
yields a difference equation for
$F_{VV}(z;p,w)=(\pi\otimes \pi)F(z;p,w)$. Noting that  $\pi\cdot\bar{\phi}_w
={\rm Ad}(D_w^{-1})\cdot \pi$, where $D_w=e_{11}+we_{22}$, we find
\beq
F_{VV}(pz;p,w)={\rm Ad}(D_w\otimes 1)(F_{VV}(z;p,w))\cdot K 
   R_{VV}(pz),\label{dif-face-sl11-image}
\eeq
where $K=(\pi\otimes\pi)q^\T=q^2 e_{11}\otimes e_{11}+q e_{11}\otimes
e_{22}+q e_{22}\otimes e_{11}+e_{22}\otimes e_{22}$ and $R_{VV}(pz)$
is given by (\ref{sl11-r}) (with $\t=\t'=1$).
(\ref{dif-face-sl11-image}) is a system of difference equations of
$q$-KZ equation type \cite{Fre92}, and can be solved with the help of
the $q$-hypergeometric series. The solution with the initial condition
(\ref{initial-sl11-face}) is given by
\bea
F_{VV}(z;p,w)&=&{}_1\phi_0(z;p,w)e_{11}\otimes e_{11}
    +e_{22}\otimes e_{22}\no\\
& & +f_{11}(z;p,w)e_{11}\otimes e_{22}
    +f_{22}(z;p,w)e_{22}\otimes e_{11}\no\\
& &+f_{12}(z;p,w)e_{12}\otimes e_{21}
    +f_{21}(z;p,w)e_{21}\otimes e_{12},\label{face-solution}
\eea
where
\bea
&&{}_1\phi_0(z;p,w)=\frac{(pq^{-2}z;p)_\infty}{(pq^2z;p)_\infty},\no\\
&&f_{11}(z;p,w)={}_2\phi_1\lt(
\begin{array}{c}
wq^{-2}~~q^{-2}\\
w
\end{array}
;p, pq^2z\rt),\no\\
&&f_{12}(z;p,w)=-\frac{w(q-q^{-1})}{1-w}\;{}_2\phi_1\lt(
\begin{array}{c}
wq^{-2}~~pq^{-2}\\
pw
\end{array}
;p, pq^2z\rt),\no\\
&&f_{21}(z;p,w)=\frac{zpw^{-1}(q-q^{-1})}{1-pw^{-1}}\;{}_2\phi_1\lt(
\begin{array}{c}
pw^{-1}q^{-2}~~pq^{-2}\\
p^2w^{-1}
\end{array}
;p, pq^2z\rt),\no\\
&&f_{22}(z;p,w)={}_2\phi_1\lt(
\begin{array}{c}
pw^{-1}q^{-2}~~q^{-2}\\
pw^{-1}
\end{array}
;p, pq^2z\rt).
\eea
Here
\bea
&&{}_2\phi_1\lt(
\begin{array}{c}
q^a~~q^b\\
q^c
\end{array}
;p, x\rt)=\sum_{n=0}^\infty\frac{q^a;p)_n(q^b;p)_n}
  {(p;p)_n(q^c;p)_n}x^n,\no\\
&&(a;p)_n=\prod_{k=0}^{n-1}(1-ap^k),~~~~(a;p)_0=1.
\eea

\subsection{Elliptic Quantum Supergroups of Vertex Type}

As we mentioned before, a given Kac-Moody superalgebras $\G$ allows
many inequivalent
simple root systems. By means of the ``extended" Weyl transformation
method introduced in \cite{Lei85},
one can transform from one simple root system to another inequivalent
one \cite{Fra89}. 
For $\G=\widehat{sl(n|n)}$, there exists a simple root system in which
all simple roots are odd (or fermionic). This system can be constructed
from the distinguished
simple root system by using the ``extended"
Weyl operation repeatedly. We find the following simple
roots, all of which are odd (or fermionic)
\bea
&&\a_0=\d-\ve_1+\d_n\,,\no\\
&&\a_{2j}=\d_j-\ve_{j+1}\,,~~~~~j=1,2,\cdots, n-1,\no\\
&&\a_{2i-1}=\ve_i-\d_i\,,~~~~~i=1, 2,\cdots,n\label{roots}
\eea
with $\d,~\{\ve_i\}_{i=1}^n$ and $\{\d_i\}_{i=1}^n$ satisfying
\bea
&&(\d,\d)=(\d,\ve_i)=(\d,\d_i)=0,~~~~(\ve_i,\ve_j)=\d_{ij},\no\\
&&(\d_i,\d_j)=-\d_{ij},~~~~(\ve_i,\d_j)=0.
\eea
Such a simple root system is usually called non-standard. 
It seems to us that $\widehat{sl(n|n)}$ is the only non-twisted affine
superalgebra which has a non-standard system of simple
roots, all of which are fermionic.

As will be shown below, 
for $\G=\widehat{sl(n|n)}$ with the above fermionic simple roots,
one can construct a different type of twistor. Following Jimbo
et al \cite{Jim97}, we say this twistor is of {\it vertex type}. 

Let us write $h_i=\a_i~(i=0,1,\cdots,2n-1)$ with $\a_i$ given by
(\ref{roots}). We extend the Cartan
subalgebra of $\widehat{sl(n|n)}$
by adding to it the element $h_{\rm ex}=\sum_{i=1}^n
(\ve_i+\d_i)$. A basis of the extended Cartan subalgebra is
$\{h_{\rm ex}, h_0,h_1,\cdots,h_{2n-1},d\}$. 
Denote by $\{h^{\rm ex},h^0,h^1,\cdots,h^{2n-1},c\}$ the dual basis. We have
\bea
&&h^{\rm ex}=\frac{1}{2n}\sum_{i=1}^n(\ve_i-\d_i),\no\\
&&h^{2k}=d+\sum_{i=1}^k(\ve_i-\d_i)-\frac{k}{n}\sum_{i=1}^n
   (\ve_i-\d_i),\no\\
&&h^{2k+1}=d+\sum_{i=1}^{k+1}\ve_i-\sum_{i=1}^k\d_i-\frac{2k+1}{2n}\sum_{i=1}^n
   (\ve_i-\d_i),
\eea
where $k=0,1,\cdots,n-1$. The canonical element $\T$ in the extended
Cartan subalgebra reads
\beq
\T=h_{\rm ex}\otimes h^{\rm ex}+\sum_{i=0}^{2n-1}(h_i\otimes h^i)
  +d\otimes c.
\eeq
Let $\tau$ be the diagram automorphism of $U_q[\widehat{sl(n|n)}]$ such that
\beq
\tau(e_i)=e_{i+1~{\rm mod}~2n},~~~~\tau(f_i)=f_{i+1~{\rm mod}~2n},~~~~
\tau(h_i)=h_{i+1~{\rm mod}~2n}.
\eeq
Obviously, the automorphism $\tau$ is non-graded since it preserves the grading
of the generators and moreover, $\tau^{2n}=1$. Then we can show
\bea
&&\tau(h_{\rm ex})=-h_{\rm ex}+\xi c,~~~~
  \tau(c)=c,~~~~\tau(h^{\rm ex})=-h^{\rm ex}+\frac{1}{2n}c,\no\\
&&\tau(h^{2k})=h^{2k+1~{\rm mod}2n}+\frac{\xi}{2n}\sum_{i=1}^n
  (\ve_i-\d_i)-\frac{\xi+n-2k-1}{2n}c,\no\\
&&\tau(h^{2k+1})=h^{2k+2~{\rm mod}2n}+\frac{\xi}{2n}\sum_{i=1}^n
  (\ve_i-\d_i)-\frac{n-2k-1}{2n}c,
\eea
where $k=0,1,\cdots,n-1$ and $\xi$ is an arbitrary constant.
Introduce element
\beq
\tilde{\rho}=\sum_{i=0}^{2n-1}h^i+\xi nh^{\rm ex},
\eeq
which gives the principal gradation
\beq
[\tilde{\rho},e_i]=e_i,~~~~[\tilde{\rho},f_i]=f_i,~~~~i=0,1,\cdots,2n-1.
\eeq
It is easily shown that
\beq
\tau(\tilde{\rho})=\tilde{\rho},~~~~~(\tau\otimes\tau)\T=\T.
\eeq
Notice also that
\bea
&&(\tau\otimes\tau)\cdot\D=\D\cdot\tau,\no\\
&&(\tau\otimes\tau)\R=\R.
\eea
Here the second relation is deduced from the uniqueness of the universal
R-matrix of $U_q[\widehat{sl(n|n)}]$. It can be shown that
\beq
\sum_{k=1}^{2n}(\tau^{k}\otimes 1)\T=\trho\otimes c+c\otimes\trho
  -\frac{2(n^2-1)-3\xi}{6}c\otimes c.
\eeq
Therefore, if we set
\beq
\tT=\frac{1}{2n}\lt(\trho\otimes c+c\otimes\trho
  -\frac{2(n^2-1)-3\xi}{6}c\otimes c\rt),
\eeq
then we have
\beq
\sum_{k=1}^{2n}(\tau^k\otimes 1)(\T-\tT)=0.\label{tau-t-t}
\eeq

Introduce an automorphism
\beq
\tphi_r=\tau\cdot{\rm Ad}\lt(q^{\frac{r+c}{n}\trho}\rt),
\eeq
which depends on a parameter $r\in\C$.
Then the $2n$-fold product
\beq
\prod_{2n\geq k\geq 1}^{\leftarrow}(\tphi_r^k\otimes 1)\lt(
  q^\tT\R\rt)^{-1}
\eeq
is a formal power series in $p^{\frac{1}{2n}}$ where $p=q^{2r}$.
Moreover, it has leading term 1 thanks to the relation (\ref{tau-t-t}).
Following Jimbo et al \cite{Jim97}, we define the vertex type twistor
\beq
E(r)=\lim_{N\rightarrow\infty}\prod_{2nN\geq k\geq 1}^{\leftarrow}
  \lt(\tphi_r^k\otimes 1\rt)
  \lt(q^\tT\R\rt)^{-1}.\label{twistor-e}
\eeq
Then one can show that
$E(r)$ satisfies the graded shifted cocycle condition
\beq
E_{12}(r)(\D\otimes 1)E(r)=E_{23}(r+c^{(1)})(1\otimes\D)E(r).
\eeq
Moreover, $E(r)$ obeys the co-unit property
\beq
(\e\otimes 1)E(r)=(1\otimes\e)E(r)=1.
\eeq
We have

\begin{Definition} {\bf (Vertex type elliptic quantum supergroup)}: We define
elliptic quantum supergroup ${\cal A}_{q,p}[\widehat{sl(n|n)}]$ 
of vertex type to be
the quasi-triangular quasi-Hopf superalgebra
$(U_q[\widehat{sl(n|n)}],\D_r,\e,\Phi(r),\R(r))$ 
together with the graded algebra
anti-homomorphism $S$ defined by (\ref{e-s}) 
and $\a_r=m\cdot(S\otimes 1)E(r)^{-1}$,
$\b_r=m\cdot(1\otimes S)E(r)$. Here $\e$ is defined by (\ref{e-s}),
and
\bea
&&\D_r(a)=E(r)\D(a)E(r)^{-1},~~~\forall a\in U_q(\G),\no\\
&&\R(r)=E(r)^T\R E(r)^{-1},\no\\
&&\Phi(r)=E_{23}(r+c^{(1)})E_{23}(r)^{-1}.
\eea
\end{Definition}

Similar to the face type case, introduce a formal parameter $\zeta$
(or spectral parameter) into $\R$ and $E(r)$ by the formulae
\bea
&&\tilde{\cal R}(\zeta)={\rm Ad}(\zeta^\trho\otimes 1)\R,\no\\
&&E(\zeta,r)={\rm Ad}(\zeta^\trho\otimes 1)E(r),\no\\
&&\tilde{\cal R}(\zeta,r)={\rm Ad}(\zeta^\trho\otimes 1)\R(r)=E(\zeta^{-1},r)^T
   \tilde{\cal R}(\zeta)E(\zeta,r)^{-1}.
\eea
Then it can be shown from the definition of $E(r)$
that $E(\zeta,r)$ satisfies the difference equation
\bea
&&E(p^\frac{1}{2n}q^{\frac{1}{n}c^{(1)}}\zeta,r)=(\tau\otimes 1)^{-1}
 (E(\zeta,r))\cdot q^\tT \tilde{\cal R}
 (p^\frac{1}{2n}q^{\frac{1}{n}c^{(1)}}\zeta),\label{dif-vertex1}\\
&&E(0,r)=1.\label{initial-vertex}
\eea
The initial condition follows from (\ref{tau-t-t}) and the fact that
we are working in the principal gradation.
(\ref{dif-vertex1}) implies that
\beq
E\lt(\lt(p^\frac{1}{2n}q^{\frac{1}{n}c^{(1)}}\rt)^{2n}\zeta,r\rt)
 =E(\zeta,r))\cdot \prod_{2n-1\geq k\geq 0}^{\leftarrow}
 \;q^\tT\; (\tau\otimes 1)^{2n-k} \tilde{\cal R}\lt(\lt(
 p^\frac{1}{2n}q^{\frac{1}{n}c^{(1)}}\rt)^{2n-k}\zeta\rt).\label{dif-vertex2}
\eeq

Some remarks are in order. In non-super case \cite{Jim97},
$\pi$ and $\tau$ are commutable in the sense that
$\pi\cdot\tau={\rm Ad}(h)\cdot\pi$ with $h$ obeying $hv_i=v_{i+1~{\rm
mod}~m}$, where $\{v_i\}$ are basis of  
the vector module $V=\C^m=\C v_1\oplus\cdots\oplus\C c_m$ of
${\cal A}_{q,p}(\hat{sl}_m)$ and
$\tau$ is the cyclic diagram automorphism of $\hat{sl}_m$. In
the super (or $\Z2$ graded) case,  
however, $\pi$ and $\tau$ are not ``commutable"
in the above sense. This is
because $\tau$ is  grading-preserving  while the
$2n$-dimensional defining representation space
$V=\C^{n|n}=\C v_1\oplus\cdots\oplus\C v_{2n}$ is graded.
So to compute the image, one has 
to work out the action of $\tau$ at the universal level and then
apply the representation $\pi$. Therefore, the knowledge of
the universal R-matrix in its explicit form is required. This
makes the image  computation of the twistor more involved
in the supersymmetric case. 

As an example, consider the simplest case of elliptic
quantum affine superalgebra ${\cal A}_{q,p}[\widehat{sl(1|1)}]$.
Let us calculate the image in the two-dimensional representation $(\pi,
V)$, ~$V=\C^{1|1}$. 
As remarked above, we have to work at the universal level
first and then apply the representation.  We have
\begin{Lemma}\label{L2}: In the principal gradation,
the action of $\tau$ on the Drinfeld generators
is represented on $V$ by
\bea
&&\tau(X^+_n)=(-1)^nz^{2n+1}q^{-n}e_{12},~~~~
  \tau(X^-_n)=(-1)^{n+1}z^{2n-1}q^{-n}e_{21},\no\\
&&\tau(H_n)=(-1)^{n+1}z^{2n}\frac{[n]_q}{n}(e_{11}+e_{22}),\no\\
&&\tau(H^{\rm ex}_n)=(-1)^{n+1}z^{2n}\frac{[2n]_q}{n}\lt(q^{-n}e_{11}
  +\frac{q-q^{-1}}{2}[n]_q(e_{11}+e_{22})\rt).
  \label{tau-rep}
\eea
\end{Lemma}
\vskip.1in
Applying $\pi\otimes
\pi$ to the both side of (\ref{dif-vertex2}) and writing
$E_{VV}(\zeta;p)\equiv (\pi\otimes\pi)E(\zeta,r)$, where
$p=q^{2r}$, we get 
\beq
E_{VV}(p\zeta;p)=E_{VV}(\zeta;p)\cdot(\pi\otimes \pi)\lt((\tau\otimes 1)
   \tilde{\R}(p^\frac{1}{2}\zeta)\rt)\cdot\tilde{\R}_{VV}(p\zeta),\label{evv}
\eeq
where $\tilde{R}_{VV}(\zeta)=(\pi\otimes\pi)\tilde{R}(\zeta)$.
In view of (\ref{tau-rep}) and the explicit formula (\ref{sl11-R}) of the
universal R-matrix, (\ref{evv}) is a system of eight difference equations.

We can also proceed directly. We have, with the help of lemma \ref{L2},
\bea
&&(\pi\otimes\pi)(\tau^{2k}\otimes 1)\lt({\rm Ad}(p^k\zeta)^\trho\otimes
    1\rt)\R^{-1}q^{-\tT}=K\cdot \bar{E}_{2k},\no\\
&&(\pi\otimes\pi)(\tau^{2k-1}\otimes 1)\lt({\rm Ad}(p^{k-\frac{1}{2}}
    \zeta)^\trho\otimes
    1\rt)\R^{-1}q^{-\tT}=\rho_{2k-1}\cdot
    K^{-1}\cdot \bar{E}_{2k-1},
\eea
where $K=(\pi\otimes \pi)q^\T$ and
\bea
\rho_{2k-1}&=&\frac{(1+q^2p^{2k-1}\zeta^2)(1+q^{-2}p^{2k-1}
    \zeta^2)}{(1+p^{2k-1}\zeta^2)^2},\no\\
\bar{E}_{2k}&=& \frac{1}{1-q^2p^{2k}\zeta^2}\lt((1-q^{-2}p^{2k}\zeta^2)
   e_{11}\otimes e_{11}+(1-q^2p^{2k}\zeta^2)e_{22}\otimes e_{22}\rt.\no\\
& &   +(1-p^{2k}\zeta^2)e_{11}\otimes e_{22}\
  +(1-p^{2k}\zeta^2)e_{22}\otimes e_{11}\no\\
& &\lt.   -(q-q^{-1})p^k\zeta
   e_{12}\otimes e_{21}
  +(q-q^{-1})p^k\zeta
   e_{21}\otimes e_{12}\rt), \label{sl11-E2k}\\
\bar{E}_{2k-1}&=& \frac{1}{1+q^{-2}p^{2k-1}\zeta^2}\lt((1+q^2p^{2k-1}\zeta^2)
   e_{11}\otimes e_{11}+(1+q^{-2}p^{2k-1}\zeta^2)e_{22}\otimes e_{22}
   \rt.\no\\
& &   +(1+p^{2k-1}\zeta^2)
   e_{11}\otimes e_{22}
  +(1+p^{2k-1}\zeta^2)e_{22}\otimes e_{11}\no\\
& &\lt.   +(q-q^{-1})p^{k-\frac{1}{2}}\zeta e_{12}\otimes e_{21}
   -(q-q^{-1})p^{k-\frac{1}{2}}\zeta
   e_{21}\otimes e_{12}\rt). \label{sl11-E2k-1}
\eea

Then
\beq
E_{VV}(\zeta;p)=\prod_{k\geq 1}^{\leftarrow}
 \rho_{2k-1}K\bar{E}_{2k}K^{-1}\bar{E}_{2k-1}
 =\rho(\zeta;p)\lt(E^1_{VV}(\zeta;p)
 +E^2_{VV}(\zeta;p)\rt),
\eeq
where
\bea
\rho(\zeta;p)&=&\frac{(-pq^2\zeta^2;p^2)_\infty}{(pq\zeta;p)_\infty
  (-pq\zeta;p)_\infty},\label{rho}\\
E^1_{VV}(\zeta;p)&=&\prod_{k\geq 1}^{\leftarrow}\frac{1}{(1+p^{2k-1}
  \zeta^2)^2}\lt((1-q^{-2}p^{2k}\zeta^2)(1+q^2p^{2k-1}\zeta^2)
  e_{11}\otimes e_{11}\rt.\no\\
& &+(1-q^2p^{2k}\zeta^2)(1+q^{-2}p^{2k-1}\zeta^2)e_{22}\otimes e_{22}
  \no\\
& &+(q-q^{-1})p^{k-\frac{1}{2}}\zeta(1-q^{-2}p^{2k}\zeta^2)e_{12}\otimes
  e_{12}\no\\
& &\lt.  -(q-q^{-1})p^{k-\frac{1}{2}}\zeta(1-q^2p^{2k}\zeta^2)
  e_{21}\otimes e_{21}\rt),\label{E1}\\
E^2_{VV}(\zeta;p)&=&\prod_{k\geq 1}^{\leftarrow}\frac{1}{1+p^{2k-1}
  \zeta^2}\lt((1-p^{2k}\zeta^2)e_{11}\otimes e_{22}
  +(1-p^{2k}\zeta^2)e_{22}\otimes e_{11}\rt.\no\\
& &\lt.-(q-q^{-1})p^k\zeta e_{12}\otimes e_{21}
  +(q-q^{-1})p^k\zeta e_{21}\otimes e_{12}\rt).\label{E2}
\eea

The infinite product in $E^2_{VV}(\zeta;p)$ can be calculated
directly and we find
\beq
E^2_{VV}(\zeta;p)=b_E(\zeta)(e_{11}\otimes e_{22}+e_{22}\otimes e_{11})
 +c_E(\zeta)(e_{12}\otimes e_{21}-e_{21}\otimes e_{12}),\label{sol-E2}
\eeq
where,
\beq
b_E(\zeta)\pm c_E(\zeta)=\frac{(pq^{\pm 1}\zeta;p)_\infty
  (-pq^{\mp 1}\zeta;p)_\infty}{(-p\zeta^2;p^2)_\infty}.
\eeq
As for $E^1_{VV}(\zeta;p)$, it can be written as
\bea
E^1_{VV}(\zeta;p)&=&X_{11}(\zeta;p) e_{11}\otimes e_{11}
  +X_{22}(\zeta;p) e_{22}\otimes e_{22}\no\\
& &  +X_{12}(\zeta;p) e_{12}\otimes e_{12}
  +X_{21}(\zeta;p) e_{21}\otimes e_{21},\label{E1=X}
\eea
where $X_{ij}(\zeta;p)$ are the solution to the
following system of four difference equations
\bea
&&X_{11}(p\zeta;p)=\frac{1}{1-q^{-2}p^2\zeta^2}\lt((1+q^{-2}p\zeta^2)
    X_{11}(\zeta;p)-p^\frac{1}{2}\zeta
    (q-q^{-1})X_{12}(\zeta;p)\rt),\no\\
&&X_{12}(p\zeta;p)=\frac{1}{1-q^{2}p^2\zeta^2}\lt(-p^\frac{1}{2}
    \zeta (q-q^{-1})X_{11}(\zeta;p)+(1+q^{2}p\zeta^2)
    X_{12}(\zeta;p)\rt),\no\\
&&X_{21}(p\zeta;p)=\frac{1}{1-q^{-2}p^2\zeta^2}\lt(p^\frac{1}{2}
    \zeta (q-q^{-1})X_{22}(\zeta;p)+(1+q^{-2}p\zeta^2)
    X_{21}(\zeta;p)\rt),\no\\
&&X_{22}(p\zeta;p)=\frac{1}{1-q^{2}p^2\zeta^2}\lt((1+q^{2}p\zeta^2)
    X_{22}(\zeta;p)+p^\frac{1}{2}\zeta
    (q-q^{-1})X_{21}(\zeta;p)\rt).\label{dif-E1}
\eea

\vskip.3in
\noindent {\bf Acknowledgements.}
Our interest in quasi-Hopf algebras was ignited by 
Bo-Yu Hou's one beautful lecture on Jimbo et al paper \cite{Jim97}
when Y.-Z.Z. was visiting Northwest University, Xi'an, in December 1997.
Y.-Z.Z. thanks Bo-Yu Hou for interesting him the subject and for
helpful discussions. The financial support from Australian Research 
Council through a Queen Elizabeth II Fellowship Grant for Y.-Z.Z is
also gratefully acknowledged.


\end{document}